\documentclass[12pt,reqno]{amsart}
\usepackage{amssymb,delarray}
\usepackage{graphicx}

\let\ti=\tilde
\def\bsy#1{\boldsymbol{#1}}
\def\scr#1{{\mathcal{#1}}}

\def\leq{\leqslant}
\def\geq{\geqslant}
\def\cc{{\mathbb C}}
\def\rr{{\mathbb R}}
\def\zz{{\mathbb Z}}
\def\pp{{\mathbb P}}
\def\nn{{\mathbb N}}
\def\pr{{\sf pr}}

\def\cp{\cc\pp}
\def\bs{\backslash}
\def\pr{\mathsf{pr}}

\newtheorem{thm}{Theorem}[section]
\newtheorem{lem}[thm]{Lemma}
\newtheorem{prop}[thm]{Proposition}
\newtheorem{claim}[thm]{Claim}
\newtheorem{df}[thm]{Definition}
\newtheorem{cor}[thm]{Corollary}

{\catcode`\@=11
\gdef\n@te#1#2{\leavevmode\vadjust{%
 {\setbox\z@\hbox to\z@{\strut#1}%
  \setbox\z@\hbox{\raise\dp\strutbox\box\z@}\ht\z@=\z@\dp\z@=\z@%
  #2\box\z@}}}
\gdef\leftnote#1{\n@te{\hss#1\quad}{}}
\gdef\rightnote#1{\n@te{\quad\kern-\leftskip#1\hss}{\moveright\hsize}}
\gdef\?{\FN@\qumark}
\gdef\qumark{\ifx\next"\DN@"##1"{\leftnote{\rm##1}}\else
 \DN@{\leftnote{\rm??}}\fi{\rm??}\next@}}

\begin{document}

\hyphenation{re-gu-lar}

\baselineskip=14.pt plus 2pt 
\abovedisplayskip = .7\abovedisplayskip
\belowdisplayskip = .7\belowdisplayskip

\title[Regular homotopy of Hurwitz curves]
{Regular homotopy of Hurwitz curves}
\author[D.~Auroux]{Denis~Auroux}
\address{Department of Mathematics\\
MIT\\
Cambridge MA 02139-4307\\
USA\\
{\rm and\ } Centre de Math\'ematiques\\
Ecole Polytechnique\\
91128 Palaiseau\\
France}
\email{auroux@math.mit.edu, auroux@math.polytechnique.fr}
\author[Vik.~S.~Kulikov]{Viktor S.~Kulikov}
\address{Steklov Mathematical Institute\\
Gubkina str., 8\\
119991 Moscow \\
Russia}
\email{kulikov@mi.ras.ru}
\author[V.~Shevchishin]{Vsevolod~V.~Shevchishin}
\address{Mathematisches Institut\\
 Abteilung f\"ur reine Mathematik \\
 Albert-Ludwigs-Universit\"at \\
             Eckerstrasse 1\\
             D-79104 Freiburg\\
Germany}
\email{sewa@email.mathematik.uni-freiburg.de}

\dedicatory{}
\thanks{This work was initiated during a stay of the authors at IPAM, UCLA.
The first author was supported in part by NSF grant DMS-0244844.
The second author was supported in part by
RFBR (No.~02-01-00786) and INTAS (No.~00-0269).}
\date{This version: November 2003. }
\begin{abstract}
We prove that any two irreducible cuspidal Hurwitz curves $C_0$ and $C_1$
(or more generally, curves with $A$-type singularities) in
the Hirzebruch surface $\bsy F_N$ with coinciding homology classes and sets
of singularities are regular homotopic; and \emph{symplectically} regular
homotopic if $C_0$ and $C_1$ are symplectic with respect to a compatible
symplectic form.
\end{abstract}

\maketitle
\setcounter{tocdepth}{2}


\def\st{{\sf st}}

\setcounter{section}{-1}
\section{Introduction}
In this paper, we deal with $J$-holomorphic curves in the
projective plane and Hurwitz curves (in particular, algebraic
curves) in the Hirzebruch surfaces $\bsy F_N$ which imitate the
behavior of plane algebraic curves with respect to pencils of
lines (the definition of Hurwitz curves is given in Section
\ref{sec2}). We restrict ourselves to the case when Hurwitz curves
can have only singularities of the types $A_n$ with $n \geq 0$
(i.e., which are locally given by $y^2=x^{n+1}$) and
also so-called negative nodes (see Section \ref{sec2}).

In \cite{Moi} Moishezon proved the existence of an infinite sequence $\bar
H_i\subset \bsy F_1$ of generic irreducible cuspidal Hurwitz curves
of degree $54$ with exactly 378 cusps and 756 nodes which have pairwise
distinct braid monodromy type. In particular, they are pairwise non-isotopic,
and almost all of them are not isotopic to an algebraic cuspidal curve.

The aim of this article is to prove the following statement.

\begin{thm} \label{m}
Any two irreducible cuspidal Hurwitz curves $\bar H_0$ and $\bar H_1$ in the
Hirzebruch surface $\bsy F_N$ having the same homology class and the same
numbers of cusps and nodes (or, in presence of negative nodes, differences
between numbers of positive and negative nodes) can be connected by a
regular homotopy.

Moreover, if $\bar H_0$ and $\bar H_1$ are symplectic with respect to some
form $\omega$ compatible with the canonical ruling of $\bsy F_N$, then
the regular homotopy between them can be made $\omega$-symplectic.
\end{thm}

A regular homotopy is a deformation family $\{ \bar H_t \}_{t\in [0,1]}$
which is an isotopy except for a finite number of values of $t$ at which the curve
undergoes the ``standard'' transformation of creation or cancellation of
a pair of nodes of opposite signs, see Section \ref{sec2} for the precise
definition.%
\medskip

The result remains true if the irreducible Hurwitz curves $\bar H_0$ and
$\bar H_1$ are allowed to present arbitrary singularities of type $A_n$. The
necessary and sufficient condition for the existence of a regular homotopy
between $\bar H_0$ and $\bar H_1$ then becomes that the numbers of
singularities of each type are the same (except in the case of nodes, for
which one should compare differences between numbers of positive and negative
nodes).

As a corollary of our main theorem, we obtain:

\begin{cor} \label{MainCor}
Let $C_0$ and $C_1$ be two ordinary cuspidal irreducible symplectic surfaces
in $(\cp^2,\omega)$, $\deg C_0=\deg C_1$, pseudoholomorphic with respect to
$\omega$-tamed almost-complex structures $J_0$ and $J_1$ respectively.
If $C_0$ and $C_1$ have the same numbers of cusps and nodes, then they can
be deformed into each other by a $C^1$-smooth symplectic regular homotopy
in $\cp^2$.
\end{cor}

The structure of the rest of this paper is as follows: in Section~\ref{sec1}
we give a symplectic isotopy result for curves in
Hirzebruch surfaces whose irreducible components are sections. In Section \ref{sec2}
we define Hurwitz curves and regular homotopy, and mention some of their elementary
properties. In Section \ref{sec3} we
introduce the main ingredient in the proof of Theorem \ref{m}, namely
braid monodromy factorizations of Hurwitz curves with $A$-type singularities.
The main results are then proved in Section~\ref{RegHom}, where the outline of
a more geometric alternative proof is also given.

Finally, in Section \ref{sec5}, we construct two so-called quasi-positive factorizations
of an element in the braid group $Br_4$ which give a negative answer to the
Generalized Garside Problem asking whether the natural homomorphism from the
semigroup of quasi-positive braids to the braid group is an embedding.
\medskip

{\it Acknowledgement.} The authors would like to express their gratitude to
the Institute for Pure and Applied Mathematics at UCLA for its hospitality
during the early stages of the preparation of this paper.

\section{Symplectic isotopy of sections}
\vskip-\baselineskip\label{sec1}
\vskip\baselineskip

Let $X$ be a Hirzebruch surface $\bsy F_k$, $k \geq 0$, and $E$ a rational
holomorphic curve with self-intersection $E ^2 = -k$, which is unique if $k\geq
1$. Furthermore, let $\omega$ be a K{\"a}hler form compatible with the complex
structure $J_{\bsy F_k}$. It is known (see \cite{Li-Liu} or \cite{La-McD})
that $(X, \omega)$ is \emph{symplecto}morphic to $\bsy F_1$ if $k$ is odd or to
$\bsy F_0$ otherwise, equipped with an appropriate K{\"a}hler structure.

Let $\scr{J}$ be the set of all $C^l$-smooth $\omega$-tame almost complex
structures on $X$ with a fixed sufficiently large non-integer $l$. This is a
Banach manifold. Let $\scr M_k$ be the total moduli space of pseudoholomorphic
curves $C$ on $X$ in the homology class $[E] \in \mathsf H^2(X, \zz)$.
Recall that $\scr M_k$ is defined as the quotient
\[
\scr M_k := \left\{ (u, J):
\vcenter{\vbox{\hsize= .42\hsize \small \parindent=0pt
$J\in \scr J$, $\;u: \pp^1 \to X$ is $J$-holo\-morphic
and represents $[E]$}}
\right\}\Big/ \bsy{PGL}(2,\cc),
\]
where the group $\bsy{PGL}(2,\cc)$ acts by holomorphic automorphisms on
$\pp^1$. By abuse of notation, we denote the points of $\scr M_k$ by $(u, J)$
or by $(C, J)$, emphasizing either the parameterizing map $u$ or the image
curve $C$. The genus formula for pseudoholomorphic curves (see \cite{Mi-Wh})
ensures that all curves in $\scr M_k$ are embedded. Let $\pr_k :\scr M_k \to
\scr J$ be the natural projection given by $\pr_k: (C, J) \mapsto J$.
This is a Fredholm operator of $\rr$-index $2(1-k)$.

\begin{lem}\label{lem1.1} The projection $\pr_k :\scr
M_k \to \scr J$ is
\begin{enumerate}
\item[(1)] an embedding of real codimension $2(k-1)$ if $k\geq1$; in
particular, it is an open embedding for $k=1$;
\item[(2)] a bundle with fiber $S^2$ over an open subset $\scr J_{\sf
reg} \subset \scr J$ if $k=0$.
\end{enumerate}
\end{lem}

\proof The properties of the projection $\pr_k$ at $(C, J)$ can be described
in terms of the \emph{normal sheaf} $\scr N_C$ of $C$ (see \cite{Iv-Sh-1} or
\cite{Iv-Sh-2}). Since $C$ is embedded, $\scr N_C$ reduces to a line
bundle, denoted by $N_C$, equipped with the Gromov operator $D_{C,J}
: L^{1,p}(C, N_C) \to L^p (C, N_C \otimes \Lambda^{(0,1)})$.

Since in our case $C$ is rational, the operator $D_{C,J}$ cannot have both
non-trivial kernel and cokernel, see \cite{H-L-S}, \cite{Iv-Sh-1}, or
\cite{Iv-Sh-2}. This means that $\pr_k$ is of ``maximal rank'' everywhere on
$\scr M_k$, i.e., either an immersion, or a submersion, or a local
diffeomorphism.

The global injectivity of $\pr_k$ in the case $k\geq1$ follows easily from the
fact that two distinct $J$-holomorphic curves $C$ and $C'$ must have
\emph{positive} intersection number, which would contradict the condition $[C]
\cdot [C'] = -k <0$.

\smallskip
In the case $k=0$ we first show that every non-empty fiber of $\pr_0$ is
compact. Assuming the contrary, we obtain a sequence of $J$-holomorphic curves
$C_\nu$ with no limit set in $\pr_0 ^{-1}(J)$. Applying Gromov's compactness
theorem we may assume that the sequence converges to some reducible curve $C^*
= \sum m_i C^*_i$. Here the image is understood in the sense of \emph{cycles}.
Since the intersection number of every $C^*_i$ with every curve in $\pr_0
^{-1}(J)$, say $C_1$, is non-negative, and the sum of these indices is zero,
we conclude that
$[C_i^*]\cdot[C_1]=0$ for every $i$.

Therefore each $C_i^*$ represents a homology class which is a positive
integer multiple of the primitive class $[C_1]$.
Counting the $\omega$-area of $C^*_i$'s we conclude that $C^*$ consists of a
unique non-multiple component homologous to $C_1$. Thus $(C^*, J)$ lies in
$\scr M_0$.

\smallskip
To show that all non-empty fibers of $\pr_0$ are diffeomorphic, it is
sufficient to prove that the complement to the image of $\pr_0$ has real
codimension $2$.  Using Gromov's compactness once more, we can find a
reducible $J$-holomorphic curve $C^* = \sum m_i C^*_i$ for every $J$ lying on
the (topological) boundary of the image of $\pr_0$. By the genus formula, all
$C^*_i$ must be embedded rational curves. Moreover, since $c_1(X) \cdot [C^*]
=2$ and the intersection form of $\bsy F_0= S^2 \times S^2$ is even, for
one of them, say $C^*_1$, we must have $c_1(X)\cdot [C^*_1] \leq 0$. Applying
the genus formula we obtain that $[C^*_1]^2 =:k \leq -2$. This implies that $J$
lies in the image of the projection $\pr_k : \scr M_k \to \scr J$, which is
locally a submanifold of codimension $\geq2$.
\qed

\smallskip
We shall denote the image $\pr_k(\scr M_k)$ by $\scr J_k$. For $k\geq1$ this is
a submanifold of $\scr J$ of real codimension $2(k-1)$, in general
\emph{not closed}.

\begin{lem}\label{lem1.2o}
Let $X$ be a Hirzebruch surface $\bsy F_k$, $\omega$ a K{\"a}hler form on $X$, $J \in
\scr J_k$ an almost-complex structure, $E$ the corresponding $J$-holomorphic curve with
self-intersection $-k$, and $C$ an irreducible $J$-holomorphic curve different
from $E$. Let also $F$ be a fiber of the ruling on $X$. Then
\begin{itemize}
\item[(i)] $C$ is homologous to $d\,[E] + f\,[F]$ with $d \geq 0$ and $f \geq
k\,d$;
\item[(ii)] $c_1(X) \cdot C = d(2-k) + 2f >0.$
\end{itemize}
\end{lem}

\proof
Fix a point $x$ on $X$ outside $E$. Observe that there exists a path
$h: [0,1] \to \scr J$, $t\in [0,1]$, which connects $J =: J_0$ with the
``standard'' structure $J_1 = J_\st$ through almost-complex structures
$J_t=h(t)$, such that $E$ is $J_t$-holomorphic for
every $t \in [0,1]$. This means that $h$ takes values in $\scr J_k$. Consider
the moduli space $\scr M_h := $
\[
\scr M_h(X, F, x) := \left\{ (C, t):
\vcenter{\vbox{\hsize= .6\hsize \small \parindent=0pt
  $t\in [0,1]$, $C$ is a rational $J_t$-holomor\-phic curve homologous to $F$
  which passes through $x$.
}}
\right\}
\]
together with the natural projection $\pr_h : \scr M_h \to [0,1]$ given by
$\pr_h : (C,t) \mapsto t$. The techniques of \cite{Bar-1} and \cite{Bar-2} ensure
that $\pr_h : \scr M_h \to [0,1]$ is a diffeomorphism. This implies that there
exists a $J$-holomorphic curve $C_0$ isotopic to $F$.

The homology group $H_2(X, \zz)$ is a free abelian group
with generators $E$ and $F$, so $[C]
= d[E] + f[F]$ with some integers $d$ and $f$. Since $C$, $E$, and $C_0$ are
$J$-holomorphic, the intersection indices $[C] \cdot [C_0]$ and $[C] \cdot [E]$ are
non-negative (see \cite{Mi-Wh}). This gives the inequalities in (i).  The
equality in part (ii) follows from formulas $c_1(X)\cdot E= 2-k$ and $c_1(X)\cdot F= 2$.
Finally, $d(2-k) + 2f = 2d + f + (f-kd)$ is non-negative and vanishes only if
$d=f =0$.
\qed

\begin{thm}\label{lem1.2} Let $X$ be a Hirzebruch surface $\bsy F_k$,
 $\omega_t$, $t\in [0,1]$ a smooth family of symplectic forms on $X$ taming the
 complex structure $J_{\bsy F_k}$, $E$ a holomorphic curve with
 self-intersection $-k$, $F$ a fiber of the ruling of $X$, and $C_0$ an
 immersed $\omega_0$-symplectic surface in $X$ such that every component $C_{0,i}$
 of $C_0$ is homologous to either to $[F]$ or to $[E] + f_i[F]$ with $f_i\geq1$.
 Assume also that the surface $C_0 \cup E$ has only positive transversal double
 points as singularities.  Then there exists an $\omega_t$-symplectic isotopy $C_t$
 between $C_0$ and a holomorphic curve $C_1$, such that each $C_t$ meets $E$
 transversally with positive intersection index.
\end{thm}

The case of main interest is when all $f_i =k$. Then $C_0$ consists of
sections of the line bundle $X\bs E\cong \scr O_{\pp^1}(k)$. The general case
corresponds to a collection of meromorphic sections with various numbers of
poles.

\proof It follows from the assumption of the theorem that both $C_0$ and $E$ are
$J_0$-holomorphic curves with respect to a same $\omega_0$-tame almost complex
structure $J_0$ on $X$ (see e.g. \cite{Iv-Sh-1}). Furthermore, there exists a path
$h: [0,1] \to \scr J$, $t\in [0,1]$, connecting $J_0$ with the ``standard''
structure $J_1 = J_\st$, such that each $J_t := h(t)$ is $\omega_t$-tame and
$E$ is $J_t$-holomorphic for every $t \in [0,1]$. This means that $h$ takes
values in $\scr J_k$.  We assume that $h$ is chosen generic enough.

By Lemma \ref{lem1.2o}, $f_i \geq k$ and $c_1(X) \cdot [C_{0,i}] = 2 + 2f_i -k$
is strictly positive. On each component $C_{0,i}$ we fix $p_i := c_1(X) \cdot
[C_{0,i}] -1 = 1 + 2f_i -k$ points $x_{i,1}, \ldots x_{i,p_i}$ in generic position.
Let $\bsy{x} =\{x_{i,j}\}$ be the whole collection of these points.  Consider
the moduli space $\scr M_h = \scr M_h(X,C_0, \bsy{x})$ of deformations of
$C_0$ as $J_t$-holomorphic curves which have the same \emph{constellation} as
$C_0$, i.e.\ such that
\begin{itemize}
\item $C$ has the same number of irreducible components as $C_0$;
\item each component $C_i$ of $C$ is rational and homologous to the corresponding
 component $C_{0,i}$ of $C_0$;
\item the component $C_i$ of $C$ passes through the same points $x_{i,1}, \ldots
 x_{i,p_i}$ as $C_{0,i}$ does.
\end{itemize}
Now, the moduli space is defined as
\[
\scr M_h(X, C_0,\bsy{x}) := \left\{ (C, t):
\vcenter{\vbox{\hsize= .45\hsize \small \parindent=0pt
$t\in [0,1]$, $C$ is a $J_t$-holomorphic curve with
the constellation of $C_0$
}}
\right\}.
\]
We denote by $\pr_h: \scr M_h \to [0,1]$ the natural projection given by
$\pr_h:(C, t) \mapsto t$. By abuse of notation we write $C \in \scr M_h$ meaning
that $(C, t)$ lies in $\scr M_h$ for some $t$.

The expected real dimension of $\scr M_h$ is $1$. The possibility to deform
the structures $h(t)$ arbitrarily near the fixed points $\bsy x$ ensures the
transversality property of the deformation problem. So $\scr M_h$ is a
manifold of the expected dimension. An important observation of \cite{H-L-S},
see also \cite{Bar-1} and \cite{Bar-2}, is that in this situation,
because the curves $C$ are rational, the
projection $\pr_h$ has no critical points.

So the statement of the theorem would follow from the properness of $\pr_h: \scr
M_h \to [0,1]$. Assuming the contrary, we would find a sequence $t_n$ converging
to some $t^* \in [0,1]$ and a sequence of $J_{t_n}$-holomorphic curves $C_n \in
\scr M_h$ with no accumulation points in $\scr M_h$. By Gromov's compactness
theorem, going to a subsequence we may assume that $C_n$ weakly converges to
some $J^*$-holomorphic curve $C^*$ with $J^* = h(t^*)$.  Since it is possible
to consider the behavior of the components of $C_n$ separately, we may assume
that $C_0$ and every $C_n$ are irreducible.  Let $C^* = \sum_{j=1} ^l m_j C^*_j$
be the decomposition of $C^*$ into irreducible components, $m_j$ being the
corresponding multiplicities. It follows from Lemma \ref{lem1.2o} that there
exists exactly one component, say $C^*_1$, which is homologous to $[E] + f^*
_1 [F]$, and every remaining component $C^*_j$ is in the homology class $f^*_j
[F]$. Moreover, if $[C_0] = [E] + f [F]$, then $\sum_{j=1} ^l m_j f^*_j =f$.
Applying the genus formula to $C^*_j$ we see that $f^*_j=1$ for every $j\geq2$.

Now recall that $C^*$ must pass through the $p = 1 + 2f -k$ marked points $x_j$
used in the definition of $\scr M_h$. On the other hand, the genericity of the
path $h(t)$ and of the points $x_j$ implies that one can have at most $p^*_1 := 1 + 2f^*_1 -k$ of the
marked points on $C^*_1$, and at most one such point on $C^*_j$ for $j \geq2$.
Altogether, this allows $C^*$ to pass through $1 + 2f^*_1 -k + l-1$ marked points, which is
strictly less than the needed $1 + 2f -k$ unless $l=1$.  But this
means that $C^*$ is irreducible and hence lies in $\scr M_h$, a contradiction.
Thus $\pr_h: \scr M_h \to [0,1]$ is proper.  \qed

\section{Hurwitz curves}\label{sec2}
\begin{df} \label{df1} A \emph{Hurwitz curve} of degree $m$ in the Hirzebruch
 surface $\bsy F_N$ is the image $\bar H=f(\mathcal S)\subset \bsy F_N$
 of an oriented closed real surface $\mathcal S$ by a smooth
 map $f:\mathcal S \to \bsy F_N\setminus E_N$ such that there exists
 a finite subset $Z\subset\bar H$ with the following properties:
\begin{itemize}
\item[(i)] The restriction of $f$ to $\mathcal S\setminus f^{-1}(Z)$ is an
 embedding, and for any $p\in \bar
 H \bs Z$, $\bar H$ and the fiber $F_{\pr(p)}$ of $\pr$ meet at $p$
 transversely with intersection index $+1;$
\item[(ii)] for each $p \in Z$ there is a neighborhood $U\subset \bsy F_N$ of $p$ such
 that $\bar H\cap U$ is a complex analytic curve, and the complex orientation of
 $\bar H\cap U \setminus \{p\}$ coincides with the orientation transported from $\mathcal
 S$ by $f$;
\item[(iii)] the restriction of\/ $\pr$ to $\bar H$ is a finite map of degree
 $m$.
\end{itemize}
\end{df}

For any Hurwitz curve $\bar H$ there is a unique minimal subset $Z\subset\bar H$
satisfying the conditions from Definition \ref{df1}. We denote it by $Z(\bar
H)$. We say that $\bar H$ is $\pr$-generic if $\pr_{\mid Z}:Z\to \pr(Z)$ is
one-to-one. A fiber of $\pr$ is \emph{$\bar H$-singular} if it meets $Z(\bar
H)$ and \emph{$\bar H$-regular} otherwise.

A Hurwitz curve $\bar H$ has an \emph{$A_k$-singularity} at $p \in Z(\bar H)$
if there is a neighborhood $U$ of $p$ and local analytic coordinates $x,y$ in
$U$ such that
\begin{itemize}
\item[(iv)] $\pr_{\mid U}$ is given by
$(x,y)\mapsto x$;
\item[(v)]  $\bar H\cap U$ is given by
$y^2=x^{k+1}$.
\end{itemize}
An ``$A_0$-singularity'' is in fact a smooth point where $\bar H$ becomes
tangent to the fiber of $\pr$; $A_1$ and $A_2$ singularities are ordinary
nodes and cusps, respectively. Therefore, we will say that $\bar H$ is
\emph{cuspidal}\/ if all its singularities are of type $A_k$ with
$0\leq k\leq 2$, and \emph{nodal}\/ if it has only $A_0$ and $A_1$
singularities.
We say that $\bar H$ has \emph{$A$-singularities}
if all its singularities
are of type $A_k$ with $k\geq0$.

For our purpose we need to extend the class of admissible singularities of
Hurwitz curves described in Definition \ref{df1} (ii) by allowing the
simplest non-holomorphic one.

\begin{df} A \emph{negative node} on a Hurwitz curve  $\bar H$ is a singular
 point $p \in Z(\bar H)$ such that
\begin{itemize}
\item[(ii$^-$)] {\it there is a neighborhood $U\subset \bsy F_N$ of $p$ such that
  $\bar H\cap U$ consists of two smooth branches meeting transversely at $p$
  with intersection index $-1$, and each branch of $\bar H\cap U$ meets the fiber
  $F_{\pr (p)}$ transversely at $p$ with intersection index $+1$.}
\end{itemize}
\end{df}

\begin{df} \label{def2}
 A Hurwitz curve $\bar H\subset \bsy F_N$ is called an almost-algebraic curve if
 $\bar H$ coincides with an algebraic curve $C$ over a disc $D(r)\subset \pp^1$ and
 with the union of $m$ pairwise disjoint smooth sections $H_{\infty, 1},\dots\break
 H_{\infty, m}$ of $\pr$ over $ \pp^1\setminus D(r)$.
\end{df}

\begin{df} \label{def3}

 Two Hurwitz curves $\bar H_0$ and $\bar H_1\subset \bsy F_N$ (possibly with
 negative nodes) are {\it $H$-isotopic} if there is a
 continuous isotopy $\phi_t:\bsy F_N\to \bsy F_N$, $t\in [0,1]$,
 fiber-preserving (i.e.\ $\exists \psi_t:\pp^1\to\pp^1$ such that
 $\pr\circ \phi_t=\psi_t\circ \pr$), and smooth outside the
 fibers $F_{\pr(s)}$, $s\in Z(\phi_t(\bar H_0))$ such that
\begin{itemize}
\item[(i)] $\phi_0=\mbox{Id}$;
\item[(ii)] $\phi_t(\bar H_0)$ is a Hurwitz  curve for all $t\in [0,1]$;
\item[(iii)] $\phi_1(\bar H_0)=\bar H_1$;
\item[(iv)]
$\phi_{t}(E_N)= E_N$ for all $t\in [0,1]$.
\end{itemize}
\end{df}

In the specific case of curves with $A$-singularities, we can in fact assume that $\phi_t$
is smooth everywhere.

The following theorem was proved in \cite{Kh-Ku}.

\begin{thm}{\rm (\cite{Kh-Ku})} Any $\pr$-generic Hurwitz curve $\bar
 H\subset \bsy F_N$ with $A$-singularities is $H$-isotopic to an almost-algebraic curve. If
 moreover $\bar H$ is a symplectic surface in $\bsy F_N$, then
 this isotopy can be chosen symplectic.
\end{thm}

\begin{df}\label{df-crea}%
 A \emph{creation of a pair of nodes} along a simple curve $\gamma$ is the
 transformation of a Hurwitz curve $\bar H=\bar{H}_{t^*-\tau}$ given by a homotopy $\bar H_t$,
 $t \in [t^*-\tau, t^* + \tau]$ and $0<\tau \ll 1$ with the following properties:
\begin{enumerate}
\item $\bar H_t$ is an isotopy outside a neighborhood $U$ of $\gamma$;
\item there exist real coordinates $(x,y,u,v)$ in $U$ such that the projection
 $\pr : \bsy F_N \to \pp^1$ is given by $\pr: (x,y,u,v) \to (x,y)$ and the curve
 $\gamma$ is given by $\{ u\in [-\tau, \tau], x=v=y=0 \}$;
\item for every $t \in [t^*-\tau, t^* + \tau]$ the curve $\bar H_t \cap U$ consists of
 two discs which are the graphs of the sections $s^\pm_t: (x,y) \mapsto (u,v)=
 \big(\pm(x^2 -(t-t^*)),\; \pm y\big)$.
\end{enumerate}%
The disc given by $D := \{ v=y=0, x\in [-\sqrt \tau, \sqrt \tau], x^2 -\tau \leq u \leq \tau - x^2 \}$ is called the
\emph{created Whitney disc}.

The inverse transformation is called the
\emph{cancellation of a pair of nodes} along the Whitney disc $D$.
\end{df}

\smallskip
\noindent\includegraphics[width=\hsize]{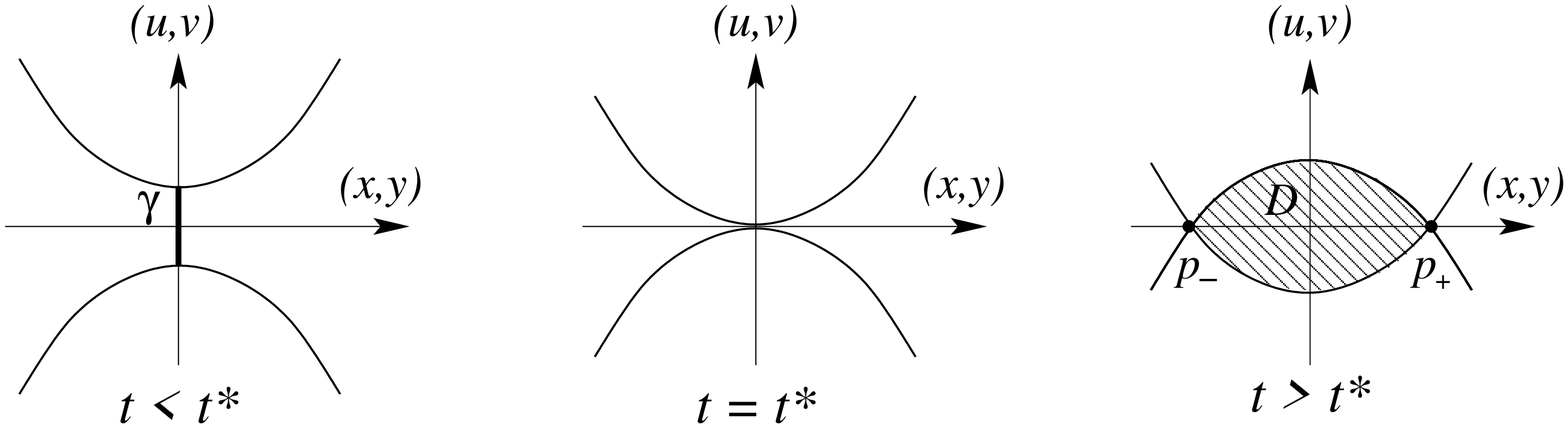}

\smallskip
The curve $\gamma$ lies in the fiber $\{x=y=0\}$ and connects two points on $\bar
H= \bar H_{t^*-\tau}$; the created nodes $p_-$ and $p_+$ have opposite
orientations: one is positive, the other negative.
The inversion of the time $t$ interchanges the creation and the cancellation
operations.

\begin{df}
 Two Hurwitz curves are \emph{regular homotopic} if one of them can be
 obtained from the other by the composition of a finite number of $H$-isotopies,
 creations and cancellations of pairs of nodes.
\end{df}

The definition is motivated by the following claim. The proof is an easy exercise.

\begin{lem}\label{gen-homo}
 Let $\phi_t: \mathcal S  \to \bsy F_N$ be a smooth homotopy of maps of
 a closed oriented real surface $\mathcal S$ to $\bsy F_N$ with the
 following properties:
\begin{enumerate}
\item for every $t$, the composition $\pr \circ \phi_t: \mathcal S  \to \pp^1$ is a
 ramified covering;
\item $\phi_t$ is an $H$-isotopy in a neighborhood of all its critical points;
\item $\phi_t$ is generic with respect to conditions $(1)$ and $(2)$.
\end{enumerate}
Then $\phi_t$ is a regular homotopy.
\end{lem}

\begin{lem}\label{create}
 Let $\bar H$ be a Hurwitz curve in $\bsy F_N$, $F$ an $\bar H$-regular fiber
 of\/ $\pr$, $\gamma \subset F$ a simple smooth curve in $F \bs \bar H$ with endpoints on
 $\bar H$, and $U$ any neighborhood of $\gamma$ in $\bsy F$. Then there exists a
 regular homotopy $\phi_t$ of $\bar H$ which creates a pair of nodes along $\gamma$
 and is constant outside~$U$.
\end{lem}

\proof It follows from the hypotheses of the lemma that $\gamma$ meets $\bar H$
transversally. So, shrinking $U$ if needed, we can find local coordinates
$(x,y,u,v)$ in $U$ satisfying the condition (2) of Definition \ref{df-crea},
such that $U \cap \bar H$ consists of two discs which are the graphs of the
mutually disjoint sections
$s^{\pm}_{t^* -\tau}: (x,y) \mapsto (u,v)= (\pm(x^2 + \tau),\pm y)$.
The result follows.
\qed

\section{Factorization semigroups}\label{sec3}
In this section we recall the notion of braid monodromy factorization
semigroups defined in \cite{Kh-Ku}.

\subsection{Semigroups over groups}\label{3.1}
A collection $(S,B,\alpha,\lambda) $, where $S$ is a semigroup, $B$ is a group, and $\alpha
:S\to B$, $\lambda :B\to \mbox{Aut}(S)$ are homomorphisms, is called {\it a semigroup
  $S$ over a group $B$} if for all $s_1,s_2\in S$
$$
s_1\cdot s_2=\lambda (\alpha (s_1))(s_2)\cdot s_1= s_2\cdot\rho (\alpha (s_2))(s_1),
$$
where $\rho (g)=\lambda(g^{-1})$. If we are given two semigroups $(S_1,B_1,\alpha_1,\lambda_1)$
and $(S_2,B_2,\alpha_2,\lambda_2)$ over, respectively, groups $B_1$ and $B_2$, we call a
pair $h=(h_S,h_B)$ of homomorphisms $h_S:S_1\to S_2$ and $h_B: B_1\to B_2$ a
\emph{homomorphism of semigroups over groups} if
\begin{itemize}
\item[(i)] $h_B\circ \alpha_{S_1}=\alpha_{S_2}\circ h_S$,
\item[(ii)] $\lambda_2(h_B(g))(h_S(s))=h_S(\lambda_1(g)(s))$
for all $s\in S_1$ and all $g\in B_1$.
\end{itemize}

The {\it factorization semigroups} defined below constitute, for
our purpose, the principal examples of semigroups over groups.

Let $\{ g_i \}_{i\in I}$ be a set of elements of a group $B$. For each $i\in I$
denote by $O_{g_i}\subset B$ the set of all the elements in $B$ conjugated to $g_i$
(the orbit of $g_i$ under the action of $B$ by inner automorphisms). Call
their union $X= \cup_{i\in I} O_{g_i}\subset B$ the {\it full set of conjugates} of $\{
g_i\}_{i\in I}$ and the pair $(B,X)$ an {\it equipped group}.

For any full set of conjugates $X$ there are two natural maps
$r=r_X: X\times X \to X$ and  $l=l_X: X\times X \to X$  defined by
$ r(a,b)=b^{-1}ab$ and  $ l(a,b)=aba^{-1}$ respectively. For each
pair of letters $a,b\in X$ denote by $R_{a,b;r}$ and $R_{a,b;l}$
the relations defined in the following way:
$$
R_{a,b;r}\quad \mbox{stands for} \quad a\cdot
b=b\cdot\,r(a,b)\quad\mbox{if $b\neq\bold 1$ and} \,\, a\cdot\bold
1=a\quad\mbox{otherwise};
$$
$$
R_{a,b;l} \quad \mbox{stands for} \quad a\cdot b=l(a,b)\cdot\,a
\quad\mbox{if $a\neq\bold 1$ and} \,\, \bold 1\cdot
b=b\quad\mbox{otherwise}.
$$
Now, put
$$
\mathcal{R}=\{ R_{a,b;r},
 R_{a,b;l} \,\vert \, (a,b)\in X\times X,\, a\neq b\, \, \mbox{if}\,\,
a\neq \bold 1 \, \, \mbox{or}\,\, b\neq \bold 1 \}
$$
and introduce a semigroup
\[
S(B,X)= \langle\, x\in X\,\,
\mid \,\, R\in \mathcal{R}\, \rangle
\]
quotient of the free semigroup of words with letters in $X$
by the relation set $\mathcal R$.
As will be seen in Section \ref{3.2}, the elements of this semigroup
represent {\it factorizations} of elements of the group $B$
with factors in $X$, up to {\it Hurwitz equivalence}.
Introduce also the {\it
product homomorphism} $\alpha_X :S(B,X)\to B$ given by $\alpha_X (x)=x$ for
each $x\in X$.

Next, we define two actions $\lambda$ and $\rho$ of the group $B$
on the set $X$:
\[
x\in X\mapsto \lambda (g) (x)=gxg^{-1}\in X.
\]
and $\rho (g)=\lambda(g^{-1})$. It is easy to see that the above
relation set $\mathcal{R}$ is preserved by the both actions;
therefore, $\rho$ and $\lambda$ define an anti-homomorphism $\rho
: B\to \mbox{Aut} (S(B,X))$ ({\it right conjugation action}) and a
homomorphism $\lambda : B\to \mbox{Aut} (S(B,X))$ ({\it left
conjugation action}). The action $ \lambda (g)$ on $S(B,X)$ is
called {\it simultaneous conjugation} by $g$.

One can easily show that  $(S(B,X),B, \alpha_X, \lambda)$ is a
semigroup over $B$. We call such semigroups the {\it factorization
semigroups over} $B$. When $B$ is a fixed group, we abbreviate
$S(B,X)$ to $S_X$. By $x_{1}\cdot\ldots\cdot x_{n}$ we denote the
element in $S_X$ defined by a word $x_{1}\dots x_{n}$.

Notice that $S: (B,X)\mapsto (S(B,X),B, \alpha_X, \lambda)$ is a functor from the category of
equipped groups to the category of semigroups over groups. In particular,
if $X\subset Y$ are two full sets of conjugates in $B$, then the identity map $id:
B\to B$ defines an embedding $id_{X,Y}:S(B,X)\to S(B,Y)$. So for each group $B$,
the semigroup $S_B=S(B,B)$ is a {\it universal factorization semigroup over}
$B$, which means that each semi-group $S_X$ over $B$ is canonically embedded
in $S_B$ by $id_{X,B}$.

Since $\alpha_X=\alpha_B \circ id_{X,B}$, we make no difference
between the product homomorphisms $\alpha_X$ and $\alpha_B$ and
denote them both simply by $\alpha$.

\begin{claim} \label{cl1} For any $s_1,\, s_2\in S(B,X)$ we have
\[
s_1\cdot s_2=s_2\cdot \rho(\alpha (s_2))(s_1)=\lambda(\alpha (s_1))(s_2)\cdot
s_1.
\]
\end{claim}

Denote by $\Sigma_m$ the symmetric group acting on the set $\{ 1, \dots,
m\}=[1,m]$ and by $(i,j)\in \Sigma_m$ the transposition exchanging
 the elements $i$ and $j\in
[1,m]$. The element $$h_g=((1,2)\cdot (1,2))^{g+1}\cdot ((2,3)\cdot(2,3))\cdot
...\cdot ((m-1,m)\cdot (m-1,m))\in S_{\Sigma_m}$$ is called a {\it Hurwitz
element of genus} $g$.

\begin{lem} \label{hur}
The Hurwitz elements $h_g$ are invariant under the conjugation action of
$\Sigma_m$ on $S_{\Sigma_m}$.
\end{lem}

\proof Put
\[
G_{h_g}=\{ \sigma\in \Sigma_m \, \mid \, \lambda(\sigma )(h_g)=
h_g\, \}.
\]
To prove Lemma \ref{hur}, it is sufficient to show that
$(i,i+1)\in G_{h_g}$ for $1\leq i\leq m-1$.

Applying Claim \ref{cl1}, and using the fact that $\lambda(\alpha(s))=\mathrm{Id}$
whenever $\alpha(s)=\bold 1$, we have
{\small
\[\displaystyle \begin{array}{l} \displaystyle
h_g=((1,2)\cdot (1,2))^{g}\cdot \prod_{j=1}^{m-1}((j,j+1)\cdot(j,j+1)) = \\
\displaystyle
((i,i+1)\cdot(i,i+1))\cdot ((1,2)\cdot (1,2))^g\cdot
\prod_{j\neq i}((j,j+1)\cdot(j,j+1))=
\\ \displaystyle
(i,i+1)\cdot \lambda ((i,i+1))\Bigl(((1,2))\cdot (1,2))^g\cdot
\prod_{j\neq i}(j,j+1)^{\cdot 2}\Bigr)
\cdot (i,i+1)= \\ \displaystyle
(i,i+1)\cdot (i,i+1) \cdot \lambda ((i,i+1))\Bigl(((1,2)\cdot (1,2))^{g}\cdot
\cdot \prod_{j\neq i}(j,j+1)^{\cdot 2}\Bigr)
= \\ \displaystyle
\lambda((i,i+1))\Bigl(((i,i+1)\cdot (i,i+1)) \cdot ((1,2)\cdot (1,2))^{g}\cdot
 \prod_{j\neq i}(j,j+1)^{\cdot 2}\Bigr)=
\\ \displaystyle \lambda((i,i+1))(h_g),
\end{array}
\]}
since $\alpha ((l,s)\cdot (l,s))=\bold 1 \in \Sigma_m$ for any transposition
$(l,s)$. \qed

\subsection{Hurwitz equivalence} \label{3.2}
As above, let $X\subset B$ be a union
of orbits of the conjugation action $\lambda$ of $B$. An ordered
set
$$\{ x_1,\dots ,x_n \,\,\mid x_i\in X \}$$
is called a {\it factorization}
(of length $n\in\nn$)
of $g= x_1\dots x_n\in B$ in $X$.
Denote by $F_{X}=\bigcup_{n\in \nn} X^n$ the set of all possible
factorizations of the elements of $B$ in $X$.
There is a natural map $\varphi : F_{X} \to S(B,X)$, given by
$$\varphi(\{ x_1,\dots ,x_n \})=x_1\cdot\dots \cdot x_n .$$
The transformations which replace in $\{ x_1,\dots ,x_n \}$
two neighboring factors $(x_i, x_{i+1})$ by
$(x_{i+1},\rho(x_{i+1})(x_i))$ or $(\lambda(x_i)(x_{i+1}),x_i)$ and
preserve the other factors are called {\it Hurwitz moves}. Two
factorizations are {\it Hurwitz equivalent} if one can be obtained
from the other by a finite sequence of Hurwitz moves.

\begin{claim} Two factorizations $x=\{ x_1,\dots ,x_n \}$ and
$x' =\{ x'_1,\dots ,x'_n \}$ are Hurwitz
equivalent if and only if $\varphi (x)=\varphi (x')$.
\end{claim}

\subsection{Semigroups over the braid group} \label{3.3}
In this subsection, $B=Br_m$ is the braid group on $m$ strings. We fix a set $\{
a_1,\dots ,a_{m-1} \} $ of {\it standard generators},
i.e., generators being subject to the relations
$$
\begin {array}{ll}
a_ia_{i+1}a_i & =a_{i+1}a_i a_{i+1} \qquad \qquad 1\leq i < m-1 ,  \\
a_ia_{k} & =a_{k}a_i  \qquad \qquad \qquad \, \, \mid i-k\mid \,
\geq 2.
\end{array}
$$

For $k\geq0$ denote by $A_k=A_k(m)$ (resp., by $A_{\bar k}=A_{\bar k}(m)$) the
full set of conjugates of $a_1^{k+1}$ (resp., of $a_1^{-k-1}$) in $Br_m$ (recall
that all the generators $a_1,\dots ,a_{m-1}$ are conjugated to each other).
Consider the factorization semigroup $S_{A_k}$ as a subsemigroup of the universal semigroup
$S_{Br_m}$ over $Br_m$.  Let $\Delta =\Delta _m $ be the so-called Garside element:
$$
\Delta= (a_1\dots a_{m-1})\dots(a_1a_2a_3)(a_1a_2)a_1.
$$
It is well-known that
$$ \Delta^2=(a_1\dots a_{m-1})^{m}$$
is the generator of the center of $Br_m$. Denote by $\delta
^2=\delta _m^2$ the element in $S_{A_0}\subset S_{Br_m}$ equal to
$$ \delta^2=(a_1\cdot \ldots \cdot a_{m-1})^{m}.$$
Put also
$$\tilde\delta_m^2=\prod_{l=m}^2 \prod^{l-1}_{k=1} z_{k,l}^2 \in S_{A_1},$$
where $z_{k,l}=(a_{l-1}\dots a_{k+1})a_k(a_{l-1}\dots
a_{k+1})^{-1}$ for $k<l$ (the notation $\prod_{l=m}^2$ means that the
product is taken over decreasing values of $l$ from left to right).
It is known (see for
example \cite{Moi-Te}) that $\alpha(\tilde\delta_m^2)=\Delta_m^2$. Moreover,
$$\tilde\delta_m^2=\prod^{m-1}_{k=1} z_{k,m}^2\cdot\tilde\delta_{m-1}^2.$$

Note that the elements $z_{k,l}^2$, $1\leq k<l\leq m$, generate the group of
{\it pure} braids $P_m= \ker \gamma_m$, where $\gamma_m:Br_m\to \Sigma_m$ is the natural
homomorphism onto the symmetric group $\Sigma_m$.

Although we will not be using that fact, it is worth mentioning that the
factorizations $\delta_m^2$ and $\tilde{\delta}_m^2$ represent the braid
monodromies of a smooth plane curve of degree $m$ and of a configuration of
$m$ lines in generic position, respectively.

\begin{lem} {\rm (\cite{Kh-Ku}) } The element $\delta ^2$ is fixed under the
conjugation action of $Br_m$ on $S_{Br_m}$, i.e., $\rho (g)(\delta
^2)=\delta ^2$ for any $g\in Br_m$.
\end{lem}

\begin{thm}  {\rm (\cite{Kh-Ku}) } \label{d2} The element $\tilde\delta_m^2$ is the only element
$s\in S_{A_1}$ such that $\alpha (s)=\Delta_m^2$.
\end{thm}

\begin{lem} \label{conj} Let $g\in A_0$ have the same image in $\Sigma_m$ as
 $z_{k_0,l_0}$, i.e., $\gamma_m(g)=\gamma_m(z_{k_0,l_0})=(k_0,l_0)$. Then there exists $p\in P_m$
 such that $g=pz_{k_0,l_0}p^{-1}$. Moreover, $p$ can be chosen so that it can
 be represented as a positive word in the alphabet $\{ z^2_{i,j}\}$.
\end{lem}

\proof The elements of $Br_m$ can be considered as diffeomorphisms of a
punctured disc $D\setminus \{x_1,\dots,x_m\}$ and, in particular, any element $g\in A_0$ can
be presented as a half-twist along a simple path $l_g$ connecting some two
points, say $x_k$ and $x_l$, in $(D\setminus \{x_1,\dots, x_m\})\cup \{ x_k,x_l\}$. It
follows from the equality $\gamma_m(g)=\gamma_m(z_{k_0,l_0})$ that $\{ x_k,x_l\} =\{
x_{k_0},x_{l_0}\}$ and the paths $l_g$ and $l_{z_{k_0,l_0}}$ are homotopic in
the closed disc $D$ as paths with fixed endpoints. Let $h_t : [0,1]\to D$, $0\leq t\leq 1$, be a
homotopy connecting $l_g$ and $l_{z_{k_0,l_0}}$, so that $h_0=l_g$
and $h_1=l_{z_{k_0,l_0}}$.

Without loss of generality, we can assume that for all $t$ the paths
$h_t([0,1])$ are simple arcs and there is a finite set $T=\{ t_1<
\dots < t_q\}$ of values $t\in (0,1)$ such that the image of $h_t$ remains
disjoint from $\{x_1,\dots,x_m\}\setminus\{x_{k_0},x_{l_0}\}$ for all
$t\not\in T$, and passes through exactly one of those points, say
$x_{j_{i}}$, for each $t_i\in T$.
Denote by $L_i$ the part of the path $h_{t_i}([0,1])$ connecting
the points $x_{k_0}$ and $x_{j_i}$.
Then it is easy to see that the paths $h_{t^+}([0,1])$, $t_i< t^+ <t_{i+1}$
are homotopic in $D\setminus \{x_1,\dots,x_m\}$ to the paths obtained from
$h_{t^-}([0,1])$, $t_{i-1}< t^-<t_{i}$ by either the full-twist along the path $L_i$
or the inverse transformation (depending on whether the map $(t,s)\mapsto
h_t(s)$ preserves or reverses orientation near $x_{j_i}$).
Since crossing each value $t_i\in T$ results in conjugation by a pure braid,
we conclude that there exists $p\in P_m$ such that $g=pz_{k_0,l_0}p^{-1}$.

To prove that $p$ can be chosen so that it can be represented as
a positive word in the alphabet $\{ z^2_{i,j}\}$, recall that
\begin{equation} \label{Delta} \Delta_m^2=\prod_{l=m}^{2} \prod^{l-1}_{k=1} z_{k,l}^2
\end{equation}
is a generator of the center of the braid group $Br_m$. Therefore a cyclic
permutation of the factors entering in (\ref{Delta}) does not change the
product. Thus, for any $z_{i_0,j_0}^2$, $\Delta_m^2$ can be presented as
$\Delta_m^2=z_{i_0,j_0}^2\Delta'_{i_0,j_0}$, where $\Delta'_{i_0,j_0}$ is a positive
word in letters $z^2_{i,j}$. This completes the proof of Lemma \ref{conj},
since
$$z_{i_0,j_0}^{-2}bz_{i_0,j_0}^2=\Delta'_{i_0,j_0}b(\Delta'_{i_0,j_0})^{-1}$$
for all $b\in Br_m$. \qed
\medskip

Consider semigroups $\mathcal A^0_k=S_{\bigcup_{i=0}^{k}A_i}$ and $\mathcal
A_k=S_{A_{\bar 1}\cup(\bigcup_{i=0}^{k}A_i)}$, $k\geq 0$. We have a natural embedding
$\mathcal A_k^0\subset \mathcal A_k$. Note that if $g\in A_1$ then $g^{-1}\in A_{\bar 1}$.
For $k\in \mathbb N\cup \{ \infty \}$, denote by
$$\overline{\mathcal A}_k = \langle\, x\in {A_{\bar 1}\cup \bigcup_{i=0}^{k}A_i}\,\, \mid \,\, R\in
\mathcal{R}_k\cup \overline{\mathcal R} \, \rangle, $$
where $\mathcal R_k$ is the set
of relations defining $\mathcal A_k$ and
 $$
\overline{\mathcal R}=\langle \, g\cdot (g^{-1})=(g^{-1})\cdot g=\bold 1\, \mid \, g\in A_1\, \rangle
$$
is the set of cancellation relations. There is a natural homomorphism of
semigroups
$$
c:\mathcal A_k\to \overline{\mathcal A}_k.
$$
We say that two elements $s_1,\, s_2\in \mathcal A_k$ are {\it weakly
 equivalent} if $c(s_1)=c(s_2)$. It is easy to see that $\overline{\mathcal
 A}_k$ can be considered as a semigroup over $Br_m$ with the natural product
homomorphism $\overline \alpha :\overline{\mathcal A}_k\to Br_m$ such that $\alpha
=\overline \alpha \circ c$.

The existence of natural embeddings of $S_{A_{\bar 1}}$ and $S_{A_j}$
($j=0,\dots,k)$ into $\mathcal A_k$
implies that any element $s\in \mathcal A_k$ can be written as
$$
s=s_{\bar 1} \cdot s_0\cdot ...\cdot s_k,
$$
where each $s_j\in S_{A_j}$ is either a product $s_j=g_{j,1}\cdot ...\cdot g_{j,d_j}$
of elements $g_{j,i}\in Br_m$ conjugated to $a_1^{j+1}$ (resp., to $a_1^{-2}$ for
$s_{\bar 1}$) or the empty word. The number $d(s_j)=d_j$ ($d_j=0$ if $s_j$ is
the empty word) is called the {\it degree} of $s_j$ in $S_{A_j}$,
$d(s)=(d_{\bar 1},d_0,\dots, d_k)$ is called the {\it multi-degree} of $s$, and
the vector $d(c(s))=(d_0,d_1-d_{\bar 1},d_2,\dots, d_k)$ is called the {\it
 multi-degree} of $c(s)\in \overline{\mathcal A}_k$.

\subsection{Braid monodromy factorizations of Hurwitz curves with $A$-singularities}
Let $\bar H$ be a $\pr$-generic Hurwitz curve with $A$-type singularities
(and possibly negative nodes) of degree $m$ in $\bsy F_N$. The classical concept of {\it braid
monodromy factorization} of $\bar H$ was given a precise definition in \cite{Kh-Ku} using
the language of factorization semigroups: a braid monodromy factorization of
$\bar H$ is an element $\text{bmf}(\bar H)\in \mathcal A_k\subset \mathcal A_{\infty}$
(for some $k\geq 0$), whose factors describe the monodromy of $\pr$ around the
various points of $Z(\bar H)$, with the property that
$\alpha(\text{bmf}(\bar H))=\Delta_m^{2N}$. If two elements $s_1,s_2\in \mathcal A_k$ are
braid monodromy factorizations of the same Hurwitz curve $\bar H$, then they
are conjugated, i.e., there exists $g\in Br_m$ such that $s_2=\lambda (g)(s_1)$. A braid
monodromy factorization $\text{bmf}(\bar H)$ of a $\pr$-generic Hurwitz
curve $\bar H$ with $A$-singularities can be written as $\text{bmf}(\bar H)=s_{\bar 1} \cdot s_0\cdot s_1\cdot ...\cdot
s_k$ with
\begin{equation}
s_{\bar 1}= \prod_{i=1}^{d_{\bar 1}}(q_{\bar 1,i} a_1^{-2}q_{\bar 1,i}^{-1})\in S_{\bar 1},
\quad\text{and}\quad
s_j=\prod_{i=1}^{d_j}(q_{j,i}a_1^{j+1}q_{j,i}^{-1})\in S_{A_j}.
\end{equation}
The index $j$ (resp., $\bar 1$) indicates the type of singularity:
the factors in $s_{\bar 1}$ are monodromies around the negative nodes,
those in $s_0$ correspond to the vertical tangency points (i.e., given by
$y^2=x$ in local coordinates), and those in $s_{j}$, $j\geq 1$, are
the singularities of type $A_j$.

The converse statement can be also proved in a straightforward
way.

\begin{thm} \label{moi} {\rm (\cite{Moi})}
For any $s\in \mathcal A_{\infty}$ such that $\alpha (s)=\Delta^{2N}$
there exists a Hurwitz curve $\bar H\subset \bsy F_N$ with
$\text{\rm bmf}(\bar H)=s$.
\end{thm}

Recall also the following statements.
\begin{thm} \label{iso} {\rm (\cite{Ku-Te}, \cite{Kh-Ku})}
Two $\pr$-generic curves with $A$-type singularities $\bar
H_1$, $\bar H_2$ in $\bsy F_N$ are $H$-isotopic if and only if
$\mbox{\rm bmf}(\bar H_1)=\mbox{\rm bmf}(\bar H_2)$. If $\bar
H_1$, $\bar H_2$ are symplectic surfaces then this
isotopy can be chosen symplectic.

\end{thm}

\begin{prop} \label{aa} {\rm (\cite{Kh-Ku})}
For any element $s \in {\mathcal A}^0_{\infty}$ such that $\alpha
(s)=\Delta _m^{2N}$ there is an almost algebraic curve $\bar
C\subset \bsy F_N$ with $\mbox{\rm bmf}(\bar C)=s$.
\end{prop}

The following theorem is a generalization of Theorem \ref{d2}:

\begin{thm} \label{d2N}
 The element $(\tilde\delta_m^2)^N$ is the only element $s\in S_{A_1}$ such that $\alpha
 (s)=\Delta_m^{2N}$.
\end{thm}

\proof As in the proof of Theorem \ref{d2}, one can show that the orbit of
$(\tilde\delta_m^2)^N$ under the conjugation action of $Br_m$ on $S_{A_1}$ consists
of the single element $(\tilde\delta_m^2)^N$ (see \cite{Kh-Ku}, Proposition
1.2).

Consider an element $s\in S_{A_1}$ such that $\alpha (s)=\Delta_m^{2N}$. By Proposition
\ref{aa}, there is an almost algebraic curve $\bar H \subset \bsy F_N$ with
$\text{bmf}(\bar H)=s$. After rescaling, we can assume that $\bar H$ is a
symplectic curve which coincides with an algebraic curve over a
disc $D\subset \pp^1$. We can change the complex structure on $\bsy
F_N$ over the complement $\pp^1\bs D$ so that it coincides with
the standard complex structure near the exceptional section $E_N$ and $\bar H$ becomes a
pseudoholomorphic curve. By Theorem \ref{lem1.2}, the curve $\bar H$ is
symplectically isotopic to an algebraic curve $\bar C$, with $\text{bmf}(\bar
C)=s$, since the isotopy preserves braid monodromy factorization (recall
that each irreducible component is a section of $\bsy F_N$).

This isotopy result implies the uniqueness of $s$ up to conjugation
by an element of $Br_m$;
hence $s$ is conjugated to $(\tilde\delta_m^2)^N$, and therefore, by the
observation made at the beginning of the proof, $s=(\tilde\delta_m^2)^N$.
Alternatively, one can also check directly that $\text{bmf}(\bar
C)=(\tilde\delta_m^2)^N$. Indeed, $\bar C$ is the union of $m$
holomorphic sections
of $\bsy F_N$ in generic position; the calculation is the same as for a
generic configuration of $m$ lines in $\cp^2$ (or can be reduced to that
case by viewing $\bsy F_N$ as a fiber sum of $N$ copies of $\bsy F_1$ which
is $\cp^2$ with one point blown up).
\qed
\medskip

\begin{lem}\label{cancel} Let $\bar H$ be a Hurwitz curve in $\bsy F_N$ and
 $p_+,\, p_- \in Z(\bar H)$ two nodal singular points, such that the fibers of
 $\bsy F_N$ through $p_\pm$ contain no other points from $Z(\bar H)$. Assume
 that there exists a simple smooth arc $\gamma \subset \pp^1$ with the following
 properties:
\begin{enumerate}
\item The endpoints of $\gamma$ are the projections $q_\pm :=
\pr(p_\pm)$ of the
 nodal points $p_\pm$;
\item The interior of $\gamma$ does not meet the branch locus
$\pr(Z(\bar H))$
 of $\bar H$;
\item Let $q_0 \in \pp^1 \bs Z(\bar H)$ be a base point, $F_0$ the
fiber of $\bsy
 F_N$ over $q_0$, $\gamma_- \subset \pp^1 \bs Z(\bar H)$ a smooth immersed path from
 $q_0$ to $q_-$, $\gamma_+ \subset \pp^1 \bs Z(\bar H)$ the path from $q_0$ to $q_+$
 obtained as the composition of $\gamma_-$ and $\gamma$, and $\mu_\pm$ the braid
 monodromies of $\bar H$ at $p_\pm$ along $\gamma_\pm$ in the braid group $Br_m(F_0)$; then $\mu_+ \cdot
 \mu_- =1 \in Br_m(F_0)$.
\end{enumerate}

Then one can cancel the pair of nodal points $(p_+, p_-)$ by a
regular homotopy $\phi_t$ constant outside any given neighborhood
$U \subset F_N$ of $\pr^{-1}(\gamma)$.

\end{lem}

\proof Set $X_\gamma := \pr^{-1}(\gamma)$, $H_\gamma := \bar H
\cap X_\gamma$, $F_q := \pr^{-1}(q)$, $H_q :=\bar H \cap F_q$, and
$\gamma^\circ := \gamma \bs \{ q_-, q_+\}$. Then $X_{\gamma}$ is
diffeomorphic to $\gamma \times F_0$ and $H_\gamma \subset
X_\gamma$ consists of $m$ smooth arcs, which we denote by
$\ti\gamma_1, \ldots \ti\gamma_m$. The assumptions imply that two
of these arcs, say $\ti\gamma_1$ and $\ti\gamma_2$ have common
endpoints $p_-$, $p_+$, at which they meet transversally, and the
remaining $m-2$ are embedded and disjoint from each other and from
$\ti\gamma_1,\, \ti\gamma_2$.

Fixing a metric on $\bsy F_N$, and considering the induced metrics
on each fiber $F_q$, we can construct families $\delta_\pm(q)$ of
arcs which depend smoothly on $q \in \gamma^\circ$ and have the
following properties:
\begin{enumerate}
\item each $\delta_\pm(q)$ is a smooth embedded arc in $F_q \bs
H_q$ with endpoints
 $\ti\gamma_1(q) := \ti\gamma_1 \cap F_q$ and $\ti\gamma_2(q) := \ti\gamma_2 \cap F_q$;
\item for $q$ close to $q_\pm$, the arc $\delta_\pm(q)$ is the
geodesic in $F_q$
 connecting $\ti\gamma_1(q)$ with $\ti\gamma_2(q)$.
\end{enumerate}

The monodromy hypothesis $(3)$ implies that for every $q \in
\gamma^\circ$ the arcs $\delta_+(q)$ and $\delta_-(q)$ are
isotopic in $F_q$ relative to their endpoints. Consequently, we
may assume that $\delta_+(q) = \delta_-(q) =: \delta(q)$ for all
$q \in \gamma^\circ$. Now we see that the union $ \sqcup_{q \in
\gamma^\circ} \delta(q) \sqcup \{ p_-, p_+\}$ forms an embedded
Whitney disk $D$ with the following properties:
\begin{enumerate}
\item $D$ is the image of a smooth embedding $f: \Delta \to \bsy
F_n$ of the disc $\Delta := \big\{ (x,u) \in \rr^2: x \in [-1,1];
-(1-x^2) \leq u \leq 1-x^2 \big\};$ \item the corner points
$p_\pm$ of $D$ are the images of the corner points
 $(\pm1, 0)$ of $\Delta$;
\item the images of the boundary arcs $I_\pm := \{ (x,u): x\in
]-1,1[, u =
 \pm(1-x^2) \}$ of $\Delta$ are the arcs $\ti\gamma_1$ and $\ti\gamma_2$.
\end{enumerate}

This makes it possible to cancel the nodes $p_+$ and $p_-$ along
the Whitney disc $D$. \qed

\begin{df} The construction described in the proof of the lemma is called the
 \emph{cancellation of the pair of nodes $p_-,\, p_+$} along the path $\gamma$.
\end{df}
We leave to the reader the proof of the fact that the cancellation
of a pair of nodes along a path $\gamma$---if it is
possible---depends only on the isotopy class of $\gamma$ relative
to its endpoints. Observe that the monodromy condition is
necessary and sufficient for the cancellability of a pair of
nodes.


It follows from Lemma \ref{cancel} that if a Hurwitz surface $\bar
H_2$ is obtained from a Hurwitz surface $\bar H_1$ by the creation
of a pair of nodes, then
$$
\text{bmf}(\bar H_2)=\text{bmf}(\bar H_1)\cdot g\cdot (g^{-1}),$$
where $g$ is some
element belonging to $A_1$. Therefore, we obtain the following statement.
\begin{claim} \label{weak}
 Two $\pr$-generic curves $\bar H_1,\, \bar H_2\subset \bsy F_N$ with $A$-singularities are regular
 homotopic in the class of Hurwitz curves with $A$-singularities and negative nodes if and
 only if
$$
c(\mbox{\rm bmf}(\bar H_1))=c(\mbox{\rm bmf}(\bar H_2)).
$$
\end{claim}

\section{Regular homotopy}
\label{RegHom}
\subsection{Weak equivalence}
\label{W-Eqval}
Let $s\in \mathcal A_k$. Consider the subgroup $(Br_m)_s$ of $Br_m$
generated by the factors of $s$,
and denote by $\gamma_s:(Br_m)_s\to \Sigma_m$ the restriction of $\gamma_m
:Br_m\to \Sigma_m$ to $(Br_m)_s$. For a given $d=(d_{\bar 1},d_0,\dots, ,d_k)$ and
$b\in Br_m$ set
$$
\mathcal A_k(d,b)=\langle \, s\in \mathcal A_k\, \mid \,
d(s)=(d_{\bar 1},d_0,\dots, ,d_k),\, \alpha(s)=b\, \rangle.
$$

\begin{thm} \label{main}  \leavevmode

\begin{itemize}
\item[(i)] For any multi-degree $d=(d_{\bar 1},d_0,\dots,d_k)$ and for any $b\in Br_m$,
the set $c(\mathcal A_k(d,b))$
consists of a finite number of elements.
\item[(ii)] Let $s', s''$ be two elements in $\mathcal A_k(d,\Delta^{2N})$,
$N\in \mathbb N$, such that $\gamma_{s'}((Br_m)_{s'})=
\gamma_{s''}((Br_m)_{s''})
=\Sigma_m$.
Then $s'$ and $s''$ are weakly equivalent.
\end{itemize}
\end{thm}

\proof Consider an element $s\in \mathcal A_k(d,b)$. It can be represented as a
product $s=s_{\bar 1}\cdot s_0\cdot ...\cdot s_k$, where $s_i=g_{i,1}\cdot ...\cdot g_{i,d_i}\in
S_{A_i}$ and all $g_{i,j}\in A_i$, i.e., all $g_{i,j}$ are conjugated to
$a_1^{i+1}$ (or $a_1^{-2}$ in the case of $g_{\bar 1,j}$; to simplify the
notation, we assume that $d_i\geq 1$ for all $i$, the adaptation to the
case $d_i=0$ is trivial and left to the reader).  By Lemma
\ref{conj}, there exists $p_{i,j}\in P_m$ such that $g_{i,j}=
p_{i,j}^{-1}z_{k_{i,j}, l_{i,j}}^{i+1}p_{i,j}$ for some $z_{k_{i,j},l_{i,j}}$.
For a fixed presentation of $s$ as the product $s=s_{\bar 1}\cdot s_0\cdot ...\cdot s_k$ with
$s_i=g_{i,1}\cdot ...\cdot g_{i,d_i}\in S_{A_i}$, put $\bar s_i=
(k_{i,1},l_{i,1})\cdot ...\cdot (k_{i,d_i},l_{i,d_i})\in S_{\Sigma_m}$.

\begin{claim} Let $s'$, $s''$ be two elements in $\mathcal A_k(d,b)$ such that for some factorizations
 $s'=s'_{\bar1}\cdot s'_0\cdot ...\cdot s'_k$ and $s''=s''_{\bar1}\cdot
 s''_0\cdot...\cdot s''_k$ the ordered collections $\overline s'=(\overline
 s'_{\bar1},\overline s'_0,\dots, \overline s'_k)$ and $\overline
 s''=(\overline s''_{\bar1},\overline s''_0,\dots, \overline s''_k)$ of
 elements from $S_{\Sigma_m}$ coincide. Then $s'$ and $s''$ are weakly equivalent.
\end{claim}

\proof As above, let $s'_i=g'_{i,1}\cdot ...\cdot g'_{i,d_i}\in S_{A_i}$, where all
$g'_{i,j}\in A_i$, and $s''_i=g''_{i,1}\cdot ...\cdot g''_{i,d_i}\in S_{A_i}$
where $g''_{i,j}\in A_i$. Since $\overline s'=\overline s''$, up to Hurwitz
equivalence we can assume that $g'_{i,j}$ and $g''_{i,j}$ are conjugated by pure
braids to powers of a same element $z_{k_{i,j},l_{i,j}}$. So, by Lemma
\ref{conj}, for all $i,j$ there exists $p_{i,j}\in P_m$ such that $g''_{i,j}= p_{i,j} g'_{i,j}
p_{i,j}^{-1}$.  Moreover, we can assume
that each $p_{i,j}=\prod_{q=1}^{M_{i,j}}z^{2}_{k_{q,i,j}l_{q,i,j}}$. Put
$M=\sum_i \sum_j M_{i,j}$.

Consider
$$
\tilde\delta_m^{-2}=\prod_{l=2}^{m} \prod^{1}_{k=l-1} z_{k,l}^{-2} \in S_{A_{\bar 1}}.
$$
Note that $c(\tilde\delta_m^{2}\cdot \tilde\delta_m^{-2})=\bold 1$ in $\overline{\mathcal A}_k$ for
any $k\geq 1$. Therefore for any $s\in \mathcal A_k$ the elements $s$ and
$s\cdot(\tilde\delta_m^{2}\cdot \tilde\delta_m^{-2})^M$ are weakly equivalent.

It follows from Claim \ref{cl1} that $\tilde\delta_m^{2}$ commutes with all
elements of $S_{Br_m}$. Therefore, after a suitable sequence of Hurwitz
moves $(\tilde\delta_m^{2})^{M_{i,j}}$ can be rewritten as
$$(\tilde\delta_m^{2})^{M_{i,j}}=
\prod_{q=1}^{M_{i,j}}z^{2}_{k_{q,i,j}l_{q,i,j}}\cdot \delta_{i,j}
$$
for some $\delta_{i,j}\in S_{A_1}$. Now
{\small
\begin{align*}
s''\cdot(\tilde\delta_m^{2}\cdot \tilde\delta_m^{-2})^M  & = \Bigl(
\prod_{i}\prod_{j=1}^{d_i}g''_{i,j}\Bigr)\cdot (\tilde\delta_m^{2}\cdot
\tilde\delta_m^{-2})^M
\\
& = \Bigl(\prod_i\prod_{j=1}^{d_i}g''_{i,j}\cdot
(\tilde\delta_m^{2})^{M_{i,j}}\Bigr)\cdot ( \tilde\delta_m^{-2})^M
\\
& =
\biggl(\prod_i\prod_{j=1}^{d_i}\Bigl((p_{i,j}g'_{i,j}p_{i,j}^{-1})\cdot
\prod_{q=1}^{M_{i,j}}z^{2}_{k_{q,i,j}l_{q,i,j}}\cdot
\delta_{i,j}\Bigr)\biggr)\cdot ( \tilde\delta_m^{-2})^M
\\
&=
\biggl(\prod_{i}\prod_{j=1}^{d_i}\Bigl(\prod_{q=1}^{M_{i,j}}
z^{2}_{k_{q,i,j}l_{q,i,j}}\cdot
g'_{i,j}\cdot \delta_{i,j}\Bigr)\biggr)\cdot
( \tilde\delta_m^{-2})^M.
\end{align*}}

 In fact, the specific contribution of positive nodes ($i=1$) can be
 treated more directly by observing that for any $j$
 the element $\tilde\delta_m^{2}$ can be written as
 $\tilde\delta_m^{2}=g_{1,j}'\cdot\delta_{1,j}$ (this follows from the
 conjugation invariance of $\tilde\delta_m^2$, see Theorem \ref{d2}), so that
 $s''_1\cdot (\tilde\delta_m^2)^{d_1}=s''_1\cdot \prod_{j=1}^{d_1}
 (g'_{1,j}\cdot \delta_{1,j})$; this allows us to decrease $M$ to
 $\sum_{i\neq 1} \sum_j M_{i,j}+d_1$.

Applying the relations $a\cdot b=b\cdot\,r(a,b)$, we can move $g'_{i,j}$
to the left and obtain that
$$s''\cdot(\tilde\delta_m^{2}\cdot \tilde\delta_m^{-2})^M=
s'\cdot \delta'\cdot(\tilde\delta_m^{-2})^M$$
for
some $\delta'\in S_{A_1}$. We have $\alpha(\delta')=\Delta^{2M}$,
since $\alpha(s')=\alpha(s'')$. It follows from Theorem \ref{d2N} that
 $\delta'=(\tilde\delta_m^{2})^{M}$. Therefore, $s'$ and $s''$ are weakly equivalent.
\qed
\medskip

Since the group $\Sigma _m$ is finite, the set of collections
$\overline s=(\overline s_{\bar 1},\overline s_0,\dots, \overline s_k)$
where each $\overline{s}_i$ is a product of $d_i$ transpositions in $S_{\Sigma_m}$
is also finite. This completes the proof of Theorem \ref{main} (i).

To prove Theorem \ref{main} (ii), as above, put
\begin{equation} \label{s1}
s'=s'_{\bar1}\cdot s'_0\cdot s'_1\cdot ...\cdot s'_k,
\end{equation}
where $s'_i=g'_{i,1}\cdot ...\cdot g'_{i,d_i}\in S_{A_i}$ and all
$g'_{i,j}\in A_i$, and similarly
\begin{equation} \label{s2}
s''=s''_{\bar1}\cdot s''_0\cdot s''_1\cdot...\cdot s''_k,
\end{equation}
where $s''_i=g''_{i,1}\cdot ...\cdot g''_{i,d_i}\in S_{A_i}$ and each $g''_{i,j}\in A_i$.
To prove (ii), it is sufficient to show that $s'$ and $s''$ have
factorizations of the form (\ref{s1}) and (\ref{s2}) such that
the collections $\overline s'=(\overline s'_{\bar1},\overline
s'_0,\dots, \overline s'_k)$ and $\overline s''=(\overline
s''_{\bar1},\overline s''_0,\dots, \overline s''_k)$ of elements of
$S_{\Sigma_m}$ coincide with each other.

Note that the image $\gamma_m(g'_{i,1})\in \Sigma_m$ (respectively, $\gamma_m(g''_{i,1})\in \Sigma_m$)
is either a transposition if $i$ is even or $\bold 1\in \Sigma_m$ if $i$ is odd. By
assumption, the elements $\gamma_m(g'_{i,j})$ (respectively, $\gamma_m(g''_{i,j})$) with
even $i$ generate the group $\Sigma_m$.

For any $a\in A_i$ and $b\in Br_m$, we have
\begin{equation} \label{even}
\gamma_m(\rho(a)(b))=\gamma_m(\lambda(a)(b))=\gamma_m(a)\gamma_m(b)\gamma_m(a)
\end{equation}
if $i$ is even and
\begin{equation} \label{odd}
\gamma_m(\rho(a)(b))=\gamma_m(\lambda(a)(b))=\gamma_m(b)
\end{equation}
if $i$ is odd.

If $i$ is odd and $r$ is even, it follows from (\ref{even}) and (\ref{odd}) that after the following composition of
Hurwitz moves
{\small
$$\begin{array}{l}
 g_{r,s}\cdot \widetilde s\cdot g_{i,j}=
\lambda(g_{r,s})(\widetilde s)\cdot \lambda(g_{r,s})(g_{i,j})\cdot g_{r,s}= \\[2pt]
\lambda(g_{r,s})(\widetilde s)\cdot \lambda(\lambda(g_{r,s})(g_{i,j}))(g_{r,s})\cdot
\lambda(g_{r,s})(g_{i,j})= \\[4pt]
 \lambda(\lambda(g_{r,s})(g_{i,j}))(g_{r,s})\cdot \rho(\lambda(\lambda(g_{r,s})
 (g_{i,j}))(g_{r,s}))(\lambda(g_{r,s})(\widetilde s)) \cdot \lambda(g_{r,s})(g_{i,j}),
\end{array}
$$}
we will have
\begin{itemize}
\item[(i)] $\gamma_m(\lambda(\lambda(g_{r,s})(g_{i,j}))(g_{r,s}))=\gamma_m(g_{r,s})$;
\item[(ii)] $\gamma_m(\rho(\lambda(\lambda(g_{r,s})(g_{i,j}))(g_{r,s}))
  (\lambda(g_{r,s})(\widetilde s_t)))=\gamma_m(\widetilde s_t)$ for each factor $\widetilde s_t$ of  $\widetilde s$;
\item[(iii)] $\gamma_m(\lambda(g_{r,s})(g_{i,j}))=
  \gamma_m(g_{r,s})\gamma_m(g_{i,j})\gamma_m(g_{r,s})$.
\end{itemize}
Therefore, since the elements $\gamma_m(g'_{r,s})$ with even $r$ generate the group $\Sigma_m$,
after a finite sequence of such Hurwitz moves we can obtain factorizations of $s'$
and $s''$ such that $\overline s_i'=\overline s_i''$ for odd $i$.
After that, we move all factors  $s_i'$ entering in $s'$ for odd
$i$ (and, similarly, all factors  $s_i''$
entering in $s''$ for odd $i$) to the left by Hurwitz moves.

Now, to prove that there exist factorizations of $s'$ and $s''$ with
$\overline s_i'=\overline s_i''$ for even $i$, note that $s'$ and $s''$
are braid monodromy factorizations of two irreducible Hurwitz
curves with $A$-singularities $f_j:\mathcal S_j \to \bar H_j\subset \bsy F_N\setminus E_N$,
$j=1,\,2$, of degree $m$ in $\bsy F_N$ (irreducibility follows from the
assumption that $\gamma_{s'}((Br_m)_{s'})=\gamma_{s''}((Br_m)_{s''})=
\Sigma_m$). Since $s'$ and $s''$ have the same multi-degree,
$\mathcal S_1$ and $\mathcal S_2$ have the same geometric genus $g$.
Applying the Hurwitz theorem
to the morphisms $\pr\circ f_j$, we have
$$
\widetilde{\gamma}_m(s')=\widetilde{\gamma}_m(s'')=h_g\in S_{\Sigma_m},
$$
where $h_g$ is the Hurwitz element of genus $g$ and $\widetilde{\gamma}_m
:\mathcal{A}_{\infty}\to S_{\Sigma_m}$ is the homomorphism induced by the natural
epimorphism $\gamma_m :Br_m\to \Sigma_m$.
To go further, we need to assign markings by integers to the various factors
in $\widetilde{\gamma}_m(s')$ and $\widetilde{\gamma}_m(s'')$.
Each factor $\gamma_m(g_{i,j}')$ of $h_g$ corresponding to $g_{i,j}'$ with even $i$
is marked by the integer $i$, and similarly for $\gamma_m(g_{i,j}'')$.
Then, to prove (ii) of Theorem \ref{main}, we must show
that $\widetilde{\gamma}_m(s')=\widetilde{\gamma}_m(s'')$ as products
with marked factors. Thus, Theorem \ref{main} (ii) follows from

\begin{claim} \label{Aur}  Two marked factorizations of $h_g$ with coinciding
 sets of marks are Hurwitz equivalent.
\end{claim}

\proof It is sufficient to show that
\[
(i,i+1)_{j_1}\cdot (i+1,i+2)_{j_2}=
(i,i+1)_{j_2}\cdot (i+1,i+2)_{j_1}
\]
For this let us perform the Hurwitz move $a\cdot b \mapsto (aba^{-1}) \cdot a$ three times.
Then we obtain
\[
\begin{array}{l}
(i,i+1)_{j_1}\cdot (i+1,i+2)_{j_2}=
(i,i+2)_{j_2}\cdot (i,i+1)_{j_1}= \\
(i+1,i+2)_{j_1}\cdot (i,i+2)_{j_2}=
(i,i+1)_{j_2}\cdot (i+1,i+2)_{j_1}
\end{array}
\]
\qed

\subsection{Proof of the main results}
Theorem \ref{m} follows from Theorems \ref{moi}, \ref{iso}, Claim
\ref{weak}, and Theorem \ref{main}. The assumption that the curves
$\bar{H}_i$ are irreducible implies that the images by $\gamma_m$ of
the factors in $\mathrm{bmf}(\bar{H}_i)$ are transpositions acting
transitively on $[1,m]$, and hence generate $\Sigma_m$ as required
in order to apply Theorem \ref{main} (ii).

To prove Corollary \ref{MainCor}, observe that,
for any cuspidal symplectic surface $C$ in $(\pp^2,
\omega)$, of degree $m$, and pseudoholomorphic with respect to some
$\omega$-tamed almost-complex structure $J$, it is possible to define
a braid monodromy factorization with respect to a generic pencil of $J$-lines
(see Section 4 of \cite{Kh-Ku}); by deforming the almost-complex structure
to the standard one,
the surface $C$ is symplectically $C^1$-smoothly isotopic in $\pp^2$ to a
symplectic surface $C'$ which becomes a Hurwitz curve in $\bsy F_1$ after
blowing up a point in $ \pp^2$. Therefore we can apply Theorem \ref{m}.

\subsection{An alternative proof}

We now describe a more geometric approach to the proof of Theorem \ref{m}.
\smallskip

{\bf Step 1.}
As above, we first reduce to the case where the combinatorics
of the branched coverings $\pr_0=\pr_{|\bar{H}_0}:\bar{H}_0\to\pp^1$ and
$\pr_1=\pr_{|\bar{H}_1}:\bar{H}_1\to\pp^1$ are the same. More precisely,
let $Z^0(\bar{H}_j)$ be the set of all points of $Z(\bar{H}_j)$ other than
nodes of any orientation, and let $\Delta_j=\pr(Z^0(\bar{H}_j))$.
Because the numbers of singular points of each type coincide, there exists
a diffeomorphism of pairs $\beta:(\pp^1,\Delta_0)\to (\pp^1,\Delta_1)$
which maps the projections of the $A_i$-singularities of $\bar{H}_0$ to
the projections of the $A_i$-singularities of $\bar{H}_1$ for all $i$.
Also choose an identification between the fibers of $\pr_0$ and $\pr_1$
above the base point in $\pp^1$.

Then the irreducibility of $\bar{H}_0$ and $\bar{H}_1$ makes it possible to
modify the diffeomorphism $\beta$ in such a way that the monodromy of
$\bar{H}_0$ along any closed loop $\gamma\in\pi_1(\pp^1-\Delta_0)$ involves
the same sheets of the $m$-fold covering $\pr_0$ as the monodromy of $\bar{H}_1$
along the closed loop $\beta(\gamma)\in\pi_1(\pp^1-\Delta_1)$. In
particular, if we consider only the $A_i$-singularities with even $i$, this
means that the symmetric group-valued monodromies
$\psi_j:\pi_1(\pp^1-\Delta_j)\to\Sigma_m$ of the $m$-fold coverings $\pr_0$
and $\pr_1$ satisfy the identity $\psi_1\circ \beta_*=\psi_0$. The argument, which
involves composing $\beta$ with a suitable sequence of elementary braids,
is essentially identical to that given in Section \ref{W-Eqval} to prove that the
ordered collections $\overline{s}'$ and $\overline{s}''$ can be assumed
to coincide, and so we will not repeat it; once again the crucial ingredient
is a Hurwitz theorem for marked factorizations in the symmetric group
(Claim \ref{Aur}).

By performing a suitable $H$-isotopy on $\bar{H}_0$ (obtained by lifting a
homotopy between $\mathrm{Id}$ and $\beta$ among diffeomorphisms of
$\pp^1$), we can restrict ourselves to the case where $\beta=\mathrm{Id}$.
\smallskip

{\bf Step 2.} Because the branching data for $\pr_0$ and $\pr_1$ agree,
there exists a Riemann surface $S$ and a covering map
$f:S\to\pp^1$, with simple branching above the points of
$\Delta_0=\Delta_1$ corresponding to $A_i$-singularities with even $i$,
such that $\bar{H}_0$ and $\bar{H}_1$ are the pushforwards via $f$ of the
graphs of smooth
sections $\sigma_0$ and $\sigma_1$ of the line bundle $L=f^*(O_{\pp^1}(N))$
over $S$. Choose local holomorphic coordinates on $\pp^1$ and
local holomorphic trivializations of $O_{\pp^1}(N)$ over small
neighborhoods of every point of $\Delta_0$, and lift them to local
holomorphic coordinates on $S$ and trivializations of $L$
near the points of $f^{-1}(\Delta_0)$.

For each $p\in\Delta_0$ corresponding to an $A_i$-singularity with $i$ even,
there exists a unique branch point of $f$, say $\bar{p}\in S$, such that
$f(\bar{p})=p$. By construction the sections of $L$ obtained by lifting
Hurwitz curves are holomorphic over a neighborhood of $\bar{p}$. A Hurwitz
curve $\bar{H}=f_*(\mathrm{graph}(\sigma))$ has an $A_i$-singularity above $p$
if and only if the jet of $\sigma$ at $\bar{p}$ (its power series expansion
in the local coordinates) has all odd degree terms
vanishing up to degree $i-1$, and has a non-zero coefficient in its degree
$i+1$ term. Define $\mathcal{C}_p\subset C^\infty(L)$ to be the space of
smooth sections of $L$ which are holomorphic over a fixed given neighborhood
of $\bar{p}$ and satisfy $\sigma^{(k)}(\bar{p})=0$ for all odd $k=1,\dots,i-1$, and
let $\mathcal{B}_p=\{\sigma\in\mathcal{C}_p,\ \sigma^{(i+1)}(\bar{p})\neq 0\}$.

For each $p\in\Delta_0$ corresponding to an $A_i$-singularity with $i$ odd,
there exist two points $p',p''\in f^{-1}(p)$ such that
$\sigma(p')=\sigma(p'')$. To make sense of this equality, recall that
$L=f^*O_{\pp^1}(N)$, so that the fibers of $L$ at $p'$ and $p''$ are canonically
isomorphic (the same is true over neighborhoods of $p'$ and $p''$).
The sections of $L$ obtained by lifting Hurwitz
curves are holomorphic over neighborhoods of $p',p''$, and
the presence of an $A_i$-singularity with $i=2r+1\ge 3$ means that
the jets of $\sigma$ at $p'$ and $p''$ coincide up to order $r$ and differ
at order $r+1$. Let $\mathcal{C}_p\subset C^\infty(L)$ be
the space of smooth sections of $L$ which are holomorphic over fixed given
neighborhoods of $p'$ and $p''$ and satisfy $\sigma^{(k)}(p')=
\sigma^{(k)}(p'')$ for all $0\le k\le r$, and
let $\mathcal{B}_p=\{\sigma\in\mathcal{C}_p,\ \sigma^{(r+1)}(p')\neq
\sigma^{(r+1)}(p'')\}$.

Finally, let $\mathcal{C}=\bigcap_{p\in\Delta_0}\mathcal{C}_p$ and
$\mathcal{B}=\bigcap_{p\in\Delta_0}\mathcal{B}_p$. Note that $\mathcal{C}$
is an infinite-dimensional complex vector space, and $\mathcal{B}$ is the
complement of a finite union of hyperplanes in $\mathcal{C}$. By construction the
sections $\sigma_0$ and $\sigma_1$ describing $\bar H_0$ and $\bar H_1$
belong to $\mathcal{B}$.

\begin{lem}
$\mathcal{B}$ is path-connected.
\end{lem}

\proof
Given $\sigma\in\mathcal{B}$ and a point $p\in \Delta_0$,
define the {\it phase} $\varphi_\sigma(p)\in \rr/2\pi\zz$ by the formulas
$\varphi_\sigma(p)=\arg(\sigma^{(i+1)}(\bar{p}))$
if $p$ is an $A_i$-singularity with $i$ even, and
$\varphi_\sigma(p)=\arg(\sigma^{(r+1)}(p'')-\sigma^{(r+1)}(p'))$
if $p$ is an $A_i$-singularity with $i=2r+1$ odd.

Given any $\sigma_0,\sigma_1\in\mathcal{B}$, choose $\theta_0$
distinct from $\pi+\varphi_{\sigma_0}(p)-\varphi_{\sigma_1}(p)$ for all
$p\in\Delta_0$, and let $\bar{\sigma}_0=e^{i\theta_0}\sigma_0$.
Considering the arc $\{e^{i\theta}\sigma_0,\ \theta\in [0,\theta_0]\}$,
it is easy to see that $\sigma_0$ and $\bar\sigma_0$ belong to the same
path-connected component of $\mathcal{B}$. Finally, since the phases of
$\bar\sigma_0$ and $\sigma_1$ are never the opposite of each other, the
line segment $\{t\bar\sigma_0+(1-t)\sigma_1,\ t\in[0,1]\}$ is entirely contained in
$\mathcal{B}$.
\qed
\medskip

{\bf Step 3.}
To conclude the argument, consider a homotopy $\sigma_t$ between $\sigma_0$
and $\sigma_1$ in $\mathcal{B}$. For generic $t$, the pushforward of the
graph of $\sigma_t$ is a Hurwitz curve $\bar{H}_t$ presenting, in addition
to the non-nodal singularities of $\bar{H}_0,\bar{H}_1$, an unspecified
number of nodes of either orientation. These correspond to pairs of points
$q,q'$ belonging to the same fiber of $f:S\to\pp^1$, distinct from
any of the previously considered points in $f^{-1}(\Delta_0)$, and such that
$\sigma_t(q)=\sigma_t(q')$. Self-intersections of higher order may also
occur sporadically in the homotopy $\{\bar{H}_t\}_{t\in [0,1]}$,
when the coincidence between
$\sigma_t(q)$ and $\sigma_t(q')$ extends to higher-order terms in the jets;
however, by
Lemma \ref{gen-homo} they can always be removed by an arbitrarily
small perturbation of $\bar{H}_t$, except in the case of creations and
cancellations of pairs of nodes. Hence we may assume that $\{\bar{H}_t\}$ is
a regular homotopy.

\section{An example of a braid having two inequivalent factorizations in $S_{A_0}$}
\label{sec5}

In this section we show that the Generalized Garside Problem has a negative
solution for $m\geq 4$, i.e., the product homomorphism $\alpha : S_{A_0(m)}\to Br_m$ is
not an embedding for $m\geq 4$.
\begin{prop}
There are two elements $s_1,\, s_2\in S_{A_0(4)}$ such that $s_1\neq s_2$, but
$\alpha (s_1)=\alpha (s_2)\in Br_4$.
\end{prop}

\proof Put $a=a_1$, $b=a_2$, $c=a_3$, where $\{ a_1,a_2,a_3\}$ is a set of
standard generators of $Br_4$.  Consider the elements
{\small
$$
\begin{array}{l}
s_1=(c^{-2}b)a(c^{-2}b)^{-1} \cdot
b \cdot (ac^3)b(ac^3)^{-1}
                  \cdot (ac^{5}b^{-1}) a (ac^{5}b^{-1})^{-1} \cdot c \cdot c, \\[3pt]
s_2=b\cdot (ac)b(ac)^{-1}\cdot (ac)b(ac)^{-1}\cdot (ac)^2b(ac)^{-2}\cdot
(ac)^2b(ac)^{-2}\cdot
(ac)^3b(ac)^{-3}.
\end{array}
$$
}

We have $s_1,s_2\in S_{A_0(4)}$.
Let us show that $\alpha (s_1)=\alpha (s_2)$.
Indeed,
{\small
\begin{align*}
\alpha(s_1)&=(c^{-2}b)a(c^{-2}b)^{-1}b(ac^3)b(ac^3)^{-1}
            (ac^{5}b^{-1})a(ac^{5}b^{-1})^{-1}c^2 \\
& = c^{-2}bab^{-1}c(cbc)ac(cbc)cb^{-1}aba^{-1}c^{-3}\\
& = c^{-2}ba(b^{-1}cbcb)ac(bcbcb^{-1})aba^{-1}c^{-3}\\
& = c^{-2}ba(cbb)ac(bbc)aba^{-1}c^{-3} \\
& = c^{-2}(bacb)^3a^{-1}c^{-3}.
\end{align*}
}
Similarly, we have
{\small
$$
\begin{array}{l}
\alpha(s_2)=b(ac)b (ac)^{-1}(ac)b(ac)^{-1}(ac)^2b (ac)^{-2}
                      (ac)^2 b(ac)^{-2} (ac)^3b(ac)^{-3}= \\[3pt]
           bacb bacb acbbacb(ac)^{-3}=(bacb)^3a^{-3}c^{-3},
\end{array}
$$}%
and to prove that $\alpha (s_1)=\alpha (s_2)\in Br_4$, it is sufficient to check that
$c(bacb)=(bacb)a$ and $a(bacb)=(bacb)c$. We have
$$
\begin{array}{l}
(bacb)a=bcaba=bcbab=cbcab=c(bacb); \\
a(bacb)=babcb=(bacb)c.
\end{array}
$$

To prove that $s_1\neq s_2$, consider the homomorphism $\widetilde{\gamma}_4 :S_{Br_4}\to
S_{\Sigma_4}$ determined by the natural epimorphism $\gamma_4 :Br_4\to \Sigma_4$.
We have
{
\small
$$ \gamma_4(a)=(1,2), \quad  \gamma_4(b)=(2,3), \quad \gamma_4(c)=(3,4).$$}
Therefore,
{\small
$$
\begin{array}{l}
\widetilde{\gamma}_4(s_1)= \\
\widetilde{\gamma}_4((c^{-2}b)a(c^{-2}b)^{-1} \cdot
b \cdot (ac^3)b(ac^3)^{-1}
                  \cdot (ac^{5}b^{-1}) a (ac^{5}b^{-1})^{-1} \cdot c \cdot c)= \\[3pt]
\widetilde{\gamma}_4(bab^{-1} \cdot
b \cdot (ac)b(ac)^{-1} \cdot (acb) a (acb)^{-1} \cdot c \cdot c)= \\[3pt]
(1,3)\cdot (2,3)\cdot (1,4)\cdot (2,4)\cdot (3,4)\cdot (3,4)
\end{array}
$$}
and
{\small
$$
\begin{array}{l}
\widetilde{\gamma}_4(s_2)= \\
\widetilde{\gamma}_4(b\cdot (ac)b(ac)^{-1}\cdot (ac)b(ac)^{-1}\cdot (ac)^2b(ac)^{-2}\cdot
(ac)^2b(ac)^{-2}\cdot
(ac)^3b(ac)^{-3})= \\[3pt]
\widetilde{\gamma}_4(b\cdot (ac)b(ac)^{-1}\cdot (ac)b(ac)^{-1}\cdot b\cdot
b\cdot
(ac)b(ac)^{-1})= \\[3pt]
(2,3)\cdot (1,4)\cdot (1,4)\cdot (2,3)\cdot (2,3)\cdot (1,4).
\end{array}
$$}
Now it is easy to see that $s_1\neq s_2$, since the groups
$(\Sigma_4)_{\widetilde{\gamma}_4(s_1)}$ and $(\Sigma_4)_{\widetilde{\gamma}_4(s_2)}$ are not isomorphic. Indeed, $(\Sigma_4)_{\widetilde{\gamma}_4(s_1)}=\Sigma_4$ and
{\small
$$(\Sigma_4)_{\widetilde{\gamma}_4(s_2)}=\{ (1,4),(2,3)\}= \mathbb Z/2\mathbb Z\times
\mathbb Z/2\mathbb Z$$} \qed

\ifx\undefined\bysame
\newcommand{\bysame}{\leavevmode\hbox to3em{\hrulefill}\,}
\fi

\end{document}